\documentclass[12pt]{amsart}
\usepackage{amsmath,amssymb,bm,color,tabu}
\usepackage{graphicx}
\usepackage{subfigure}
\usepackage{diagbox}
\usepackage{enumerate}
\newtheorem{theorem}{Theorem}[section]
\newtheorem{lemma}[theorem]{Lemma}
\newtheorem{remark}[theorem]{Remark}
\numberwithin{equation}{section}
\numberwithin{table}{section}
\numberwithin{figure}{section}
\allowdisplaybreaks[4]
\begin{document}
\def\argmin{\mathop{\rm argmin}}
\def\bK{\mathbb{K}}
\def\LTG{{L_2(\Gamma)}}
\def\O{\Omega}
\def\p{\partial}
\def\R{\mathbb{R}}
\def\G{\Gamma}
\def\LT{{L_2(\O)}}
\def\HOne{{H^1(\O)}}
\def\bn{\bm{n}}
\def\bw{\bm{w}}
\def\bt{\bm{t}}
\def\cT{\mathcal{T}}
\def\tu{\tilde u}
\def\ty{\tilde y}
\def\CPF{C_{\rm PF}}
\def\CTR{C_{\rm Tr}}
\def\COM{\CPF\max(\alpha^{-1},1)}
\def\HH{{H^{1/2}(\G)}}
\def\MS{V^{ms}_{H,h}}
\def\MSv{v^{ms}_{H,h}}
\def\MSw{w^{ms}_{H,h}}
\def\MSz{z^{ms}_{H,h}}
\def\bMSy{{\bar y}_{H,h,\rho}^{ms}}
\def\bMSu{{\bar u}_{H,h,\rho}^{ms}}
\def\bMSp{{\bar p}_{H,h,\rho}^{ms}}
\def\RMS{R_{H,h}^{ms}}
\def\cE{\mathcal{E}}
\def\IH{I_{\scriptscriptstyle H}}

\title[A New Error Analysis for Neumann Boundary Control Problems]
{A New Error Analysis for Finite Element Methods for Elliptic
 Neumann Boundary Control Problems
with Pointwise Control Constraints}
\author{Susanne C. Brenner}
\address{S.C. Brenner, Department of Mathematics and Center for
Computation \& Technology,
 Louisiana State University, Baton Rouge, LA 70803, USA}
\email{brenner@math.lsu.edu}
\author{Li-yeng Sung}
\address{L.-Y. Sung, Department of Mathematics and Center for
Computation \& Technology,
 Louisiana State University, Baton Rouge, LA 70803, USA}
\email{sung@math.lsu.edu}
\thanks{This work  was supported in part
 by the National Science Foundation under
 Grant No. DMS-22-08404}
\begin{abstract}
We present a new error analysis for finite element methods for a linear-quadratic
 elliptic optimal control problem with
Neumann boundary control and pointwise control constraints.
  It can be applied to standard finite element methods when
the coefficients in the elliptic operator are smooth and also to
multiscale finite element methods when the coefficients are rough.
\end{abstract}
\date{October 25, 2024}
\maketitle
\section{Introduction}
 Let $\O\subset \R^d$ ($d=2,3$) be a Lipschitz polyhedral domain
 with boundary $\G$,
 $y_d,f\in\LT$, $\phi_1,\phi_2\in \LTG$, $\phi_1\leq\phi_2$ on $\G$,
 $\gamma\in (0,1]$ be a positive constant,
 and
\begin{equation}\label{eq:aDef}
  a(y,z)=\int_\O A\nabla y\cdot\nabla z\,dx+\int_\O \kappa yz\,dx,
\end{equation}
 where the $d\times d$ symmetric
 coefficient matrix $A(x)$ is positive definite,
  and $\kappa$ is a nonnegative bounded measurable function such
  that $\|\kappa\|_{L_1(\O)}>0$.
 More precisely, we assume that the components of $A(x)$
 are Lebesgue measurable functions and
 there exist two positive numbers
 $\alpha\leq\beta$ such that
\begin{equation}\label{eq:RoughCoefficient}
  \alpha |\xi|^2\leq \xi^t A(x)\xi\leq \beta |\xi|^2\qquad
  \forall\,x\in \O,\;\xi\in\R^d.
\end{equation}
\par
 The model linear-quadratic Neumann boundary control problem
  (cf. \cite{Lions:1971:OC,Troltzsch:2010:OC})
 is to find
\begin{equation}\label{eq:LQ}
  (\bar y,\bar u)=\argmin_{(y,u)\in\bK}\frac12\big[\|y-y_d\|_\LT^2+
  \gamma\|u\|_\LTG^2\big],
\end{equation}
 where $(y,u)$ belongs to $\bK\subset \HOne\times \LTG$ if and only if
\begin{equation}\label{eq:PDEConstraint}
  a(y,z)=\int_\O fz\,dx+\int_\Gamma uz\,ds \qquad\forall\,z\in\HOne, \\
\end{equation}
 and
\begin{equation}\label{eq:ControlConstraint}
  \phi_1\leq u \leq \phi_2 \quad\text{on}\quad \Gamma.
\end{equation}
 Here $ds$ is the infinitesimal arc length ($d=2$) or the infinitesimal
 surface area ($d=3$).
\begin{remark}\label{rem:Notation}\rm
  We follow the standard notation for differential operators, function
  spaces and  norms that can be  found for example in
  \cite{Ciarlet:1978:FEM,ADAMS:2003:Sobolev,BScott:2008:FEM}.
\end{remark}
\begin{remark}\label{rem:BoundedAndCoercive}\rm
 It is clear that the bilinear form $a(\cdot,\cdot)$ is bounded on $\HOne$.
  From a generalized Poincar\'e-Friedrichs inequality
  (cf. \cite{Necas:2012:Direct})  we also have
\begin{equation*}
  \|v\|_\HOne^2\leq \CPF\big(|v|_\HOne^2+\|\sqrt\kappa v\|_\LT^2\big),
\end{equation*}
 which together with \eqref{eq:RoughCoefficient} implies the
 coercivity of $a(\cdot,\cdot)$:
\begin{equation}\label{eq:Coercive}
  \|v\|_\HOne^2\leq \CPF\max(\alpha^{-1},1)a(v,v)\qquad\forall\,v\in\HOne.
\end{equation}
 Therefore \eqref{eq:PDEConstraint} is uniquely solvable.
\end{remark}
\begin{remark}\label{rem:Interpretation}\rm
  We can interpret \eqref{eq:PDEConstraint} as  a steady state
  diffusion-reaction model for a chemical
  substance in a heterogeneous media with the source $f$ and the
  boundary control (in a weak sense)
\begin{equation*}
  (A\nabla y)\cdot \bn=u \qquad\text{on}\quad\Gamma,
\end{equation*}
 where $\bn$ is the unit outer normal on $\Gamma$.
\end{remark}
\par
 There is a substantial literature (cf.  the monographs
 \cite{NST:2006:Elliptic,LY:2008:Adaptive,HPUU:2009:Book}
 and the references therein)
 on the analysis of finite element methods for elliptic optimal control
 problems beginning with the
 pioneering work of Falk in \cite{Falk:1973:Control}.
 Nevertheless the analyses in the literature are not directly applicable to
 multiscale finite element methods for problems with rough coefficients,
  i.e., when \eqref{eq:RoughCoefficient} is the only assumption on the
  coefficient matrix $A(x)$.
\par
 Recently a new error analysis for distributed elliptic optimal control problems
 with pointwise
 control constraints was developed in
  \cite{BGS:2024:LOD_CC,BSung:2024:NewAnalysisArXiv}
 that is applicable
 to standard finite element methods (when $A$ is smooth)
 and also to multiscale finite element methods (when $A$ is rough).
  Our goal in this paper is to develop a similar error
 analysis for the model Neumann boundary control problem defined by
 \eqref{eq:LQ}--\eqref{eq:ControlConstraint}.
\par
 The rest of the paper is organized as follows.  We recall the properties of the
 continuous problem in Section~\ref{sec:Continuous} and introduce the
 discrete problem in Section~\ref{sec:Discrete}.  An abstract error
  analysis is given in Section~\ref{sec:Analysis},
 followed by application to problems with smooth coefficients in 
 Section~\ref{sec:Smooth} and
 application to problems with rough coefficients in
 Section~\ref{sec:Rough}.  We end with some concluding
 remarks in Section~\ref{sec:Conclusions}.  Some estimates for the
 continuous problem are given in Appendix~\ref{appendix:Bounds}.
\par
Throughout the paper we use $C$ (with or without subscript) to denote
a generic positive constant that can take different values at different occurrences.
\section{The Continuous Problem}\label{sec:Continuous}
 The convex minimization problem defined by
 \eqref{eq:LQ}--\eqref{eq:ControlConstraint}
 has a unique solution $(\bar y,\bar u)$ characterized by the
 first order optimality condition (cf. \cite{ET:1999:Convex})
\begin{equation}\label{eq:VI}
  \int_\O (\bar y-y_d)(y-\bar y)dx+\gamma\int_\Gamma
  \bar u(u-\bar u)ds\geq0   \qquad\forall\, (y,u)\in\bK.
\end{equation}
\par
 Let the adjoint state $\bar p\in\HOne$ be defined by
\begin{equation}\label{eq:AdjointState}
  a(q,\bar p)=\int_\O (\bar y-y_d)q\,dx\qquad\forall\,q\in\HOne.
\end{equation}
 We can deduce from \eqref{eq:PDEConstraint} and
 \eqref{eq:AdjointState} that
\begin{align*}
  \int_\O (\bar y-y_d)(y-\bar y)dx&=a(y-\bar y,\bar p)\\
  &=\Big(\int_\O f\bar p\, dx+\int_\Gamma u\bar p\,ds\Big)-
     \Big(\int_\O f\bar p\, dx+\int_\Gamma \bar u\bar p\,ds\Big)\\
 &=\int_\Gamma(u-\bar u)\bar p\,ds \hspace{180pt}\forall\,(y,u)\in\bK.
\end{align*}
 Consequently we have, by \eqref{eq:VI},
\begin{equation}\label{eq:NewInequality}
  \int_\Gamma (\bar p+\gamma\bar u)(u-\bar u)ds\geq 0
  \qquad \forall\,u\in U_{ad},
\end{equation}
 where
\begin{equation}\label{eq:AdmissibleControl}
  U_{ad}=\{u\in \LTG:\,\phi_1\leq u\leq \phi_2\}.
\end{equation}
\par
 It follows from \eqref{eq:NewInequality} that $\bar u$ is
 the $\LTG$ orthogonal
 projection of $-(1/\gamma)\bar p$ on
 $U_{ad}$ and hence
\begin{equation}\label{eq:ControlAndAdjointState}
  \bar u=\max\big(\phi_1,\min(\phi_2,-(1/\gamma)\bar p)\big).
\end{equation}
\par
 Let the functions $\lambda_1,\lambda_2\in \LTG$ be defined by
\begin{equation}\label{eq:Lambda1}
  \lambda_1=\begin{cases}
  \bar p+\gamma \phi_1&\qquad\text{if $\phi_1\geq -(1/\gamma)\bar p$}\\[4pt]
    0 &\qquad\text{if $\phi_1<-(1/\gamma)\bar p$}
  \end{cases} \quad =\max(\bar p+\gamma\phi_1,0),
\end{equation}
 and
\begin{equation}\label{eq:Lambda2}
  \lambda_2=\begin{cases}
    \bar p+\gamma \phi_2&\qquad\text{if $\phi_2\leq -(1/\gamma)\bar p$}\\[4pt]
    0 &\qquad\text{if $\phi_2>-(1/\gamma)\bar p$}
  \end{cases} \quad =\min(\bar p+\gamma\phi_2,0).
\end{equation}
 Then it follows from \eqref{eq:ControlAndAdjointState}--\eqref{eq:Lambda2}
 (cf. \cite[Section~2.8.4]{Troltzsch:2010:OC}) that
\begin{align}
  \bar p+\gamma\bar u&=\lambda_1+\lambda_2, \label{eq:KKT}\\
  \lambda_1\geq0\quad&\text{and}\quad \lambda_2\leq0, \label{eq:Signs}\\
  \int_\Gamma\lambda_1(\bar u-\phi_1)ds&=0=\int_\Gamma\lambda_2(\bar u-\phi_2)ds.
  \label{eq:Complemntarity}
\end{align}
%
\section{The Discrete Problem}\label{sec:Discrete}
 Let $\cT_\rho$ be a regular triangulation of $\Gamma$ with
  mesh size $\rho$ and $V_*$
 be a finite dimensional subspace of $\HOne$.  The space of piecewise constant
 functions with respect to $\cT_\rho$ is denoted by $W_\rho$, and the orthogonal
 projection operator from $\LTG$ onto $W_\rho$ is denoted by
 $Q_\rho$.
\par
 The discrete problem for \eqref{eq:LQ}--\eqref{eq:ControlConstraint} is to find
\begin{equation}\label{eq:DP}
  (\bar y_{*,\rho},\bar u_{*,\rho})=\argmin_{(y_*,u_\rho)\in \bK_{*,\rho}}
  \frac12\big[\|y_*-y_d\|_\LT^2+\gamma\|u_\rho\|_\LTG^2\big].
\end{equation}
 Here $(y_*,u_\rho)$ belongs to $\bK_{*,\rho}$ if and only if $u_\rho\in U_{ad}^\rho$, where
\begin{equation}\label{eq:DiscreteControlConstraint}
   U_{ad}^\rho=\{w_\rho\in W_\rho:\,Q_\rho \phi_1\leq w_\rho\leq Q_\rho\phi_2\},
\end{equation}
 and
\begin{equation}\label{eq:y*}
  y_*=y_*^\dag-B_* u_\rho,
\end{equation}
 where $y_*^\dag\in V_*$ satisfies
\begin{equation}\label{eq:DiscretePDEConstraint}
  a(y_*^\dag,z_*)=\int_\O fz_*dx+\int_\Gamma u_\rho z_*ds
  \qquad\forall\,z_*\in V_*.
\end{equation}
 The operator $B_*$ is a bounded linear operator from $\LTG$ to $\HOne$
  whose image is
 orthogonal to $V_*$ with respect to $a(\cdot,\cdot)$, i.e.,
\begin{equation}\label{eq:Orhtogonality}
  a(B_* v,z_*)=0 \qquad\forall\,v\in\LTG,\,z_*\in V_*.
\end{equation}
\begin{remark}\label{rem:B}\rm
  The operator $B_*=0$ for standard finite element methods.  For multiscale
  finite element methods
  $B_*$ is a correction operator for the Neumann boundary data.
\end{remark}
\begin{remark}\label{rem:AdmissibleSets}\rm
  It follows from \eqref{eq:AdmissibleControl} and
  \eqref{eq:DiscreteControlConstraint} that
  $Q_\rho$ maps $U_{ad}$ onto $U_{ad}^\rho$.
  In particular $U_{ad}^\rho$ is nonempty.
\end{remark}
\par
 The convex minimization problem defined by
 \eqref{eq:DP}--\eqref{eq:DiscretePDEConstraint}
 has a unique solution $(\bar y_{*,\rho},\bar u_{*,\rho})$ in $\bK_{*,\rho}$
 characterized by the first order optimality condition
\begin{equation}\label{eq:DiscreteFirstOrderOptimality}
  \int_\O (\bar y_{*,\rho}-y_d)(y_*-\bar y_{*,\rho})dx+
  \gamma\int_\Gamma \bar u_{*,\rho}(u_\rho-\bar u_{*,\rho})ds\geq0
  \qquad\forall\,(y_*,u_\rho)\in\bK_{*,\rho}.
\end{equation}
\par
 The discrete adjoint state $\bar p_{*,\rho}\in V_*$ is defined by
\begin{equation}\label{eq:DiscreteAdjointState}
  a(q_*,\bar p_{*,\rho})=\int_\O (\bar y_{*,\rho}-y_d)q_*dx \qquad\forall\,q_*\in V_*.
\end{equation}
\par
 The following simple stability result is useful for the error analysis of the discrete problem.
\begin{lemma}\label{lem:SimpleBdd}
   Let $V$ be a closed subspace of $\HOne$, $r\in\LT$ and $g\in\LTG$.
   If $v\in V$ is defined by
\begin{equation}\label{eq:VP}
  a(v,w)=\int_\O rw\,dx+\int_\G g w\,ds \qquad\forall\,w\in V,
\end{equation}
 then we have
\begin{align}
  \|v\|_\HOne&\leq \CPF\max(\alpha^{-1},1)\big(\|r\|_\LT
  +\CTR\|g\|_\LTG\big),\label{eq:SimpleBound}\\
  \|v\|_a&\leq \sqrt{\COM}\big(\|r\|_\LT+\CTR\|g\|_\LTG\big),
  \label{eq:SimpleBound2}
\end{align}
 where $\CTR$ is the positive constant that appears in the trace inequality
\begin{equation}\label{eq:Trace}
  \|\zeta\|_\LTG\leq \CTR \|\zeta\|_\HOne \qquad\forall\,\zeta\in\HOne.
\end{equation}
\end{lemma}
\begin{proof}
  It follows from \eqref{eq:Coercive}, \eqref{eq:VP} and \eqref{eq:Trace}
  that
\begin{align*}
  \|v\|_\HOne^2&\leq \COM a(v,v)\\
  &=\COM\Big(\int_\O rv\,dx+\int_\G gv\,ds\Big)\\
  &\leq \COM\big(\|r\|_\LT\|v\|_\LT+\|g\|_\LTG\CTR\|v\|_\HOne\big),
\end{align*}
 which implies \eqref{eq:SimpleBound}.
\par
 Similarly we have, by \eqref{eq:Coercive}, \eqref{eq:VP} and \eqref{eq:Trace},
\begin{align*}
  \|v\|_a^2&=\int_\O rv\,dx+\int_\G gv\,ds\\
  &\leq (\|r\|_\LT+\CTR\|g\|_\LTG)\|v\|_\HOne\\
  &\leq (\|r\|_\LT+\CTR\|g\|_\LTG)\sqrt{\COM}\|v\|_a,
\end{align*}
 which implies \eqref{eq:SimpleBound2}.
\end{proof}
\section{An Abstract Error Analysis}\label{sec:Analysis}
 We will derive error estimates in terms of the correction operator
 $B_*:\LTG\longrightarrow\HOne$, the  $L_2$ projection operator
  $Q_\rho:\LTG\longrightarrow W_\rho$, and the Ritz projection operator
  $R_*:\HOne\longrightarrow V_*$ defined by
\begin{equation}\label{eq:Ritz}
  a(R_*\zeta,v_*)=a(\zeta,v_*) \qquad\forall\,v_*\in V_*.
\end{equation}
\begin{theorem}\label{thm:LTwoError}
  There exists a positive constant $C_\flat$ depending only on $\alpha^{-1}$
  and $\gamma^{-1}$ such that
\begin{align*}
  &\|\bar y-\bar y_{*,\rho}\|_\LT^2+\|\bar u-\bar u_{*,\rho}\|_\LTG^2\notag\\
 &\hspace{30pt} \leq C_\flat\Big(\|\bar y-(R_*\bar y-B_*\bar u)\|_\LT^2
  +(1+\|B_*\|^2)\|\bar u-Q_\rho\bar u\|_\LTG^2
    +\|\bar p-R_*\bar p\|_\LTG^2\notag\\
    &\hspace{60pt} +\|B_*\|^2\|\bar p-R_*\bar p\|_a^2
  +\|\lambda_1-Q_\rho\lambda_1\|_\LTG^2 +\|\lambda_2-Q_\rho\lambda_2\|_\LTG^2\\
   &\hspace{90pt}
   +\|\phi_1-Q_\rho\phi_1\|_\LTG^2+\|\phi_2-Q_\rho\phi_2\|_\LTG^2\Big),\notag
\end{align*}
 where $\|B_*\|$ is the operator norm of $B_*:\LTG\longrightarrow\HOne$
 with respect to the norms $\|\cdot\|_\LTG$ and $\|\cdot\|_a$.
\end{theorem}
\begin{proof}
  From \eqref{eq:AdjointState}, \eqref{eq:y*}--\eqref{eq:Orhtogonality}
    and \eqref{eq:Ritz}, we find, for any $(y_*,u_\rho)\in\bK_{*,\rho}$,
\begin{align}\label{eq:LTwoError1}
  \int_\O (\bar y-y_d)(y_*-\bar y_{*,\rho})dx&=a(y_*-\bar y_{*,\rho},\bar p)\notag\\
  &=a(y_*^\dag-\bar y_{*,\rho}^\dag,\bar p)
  +a\big(B_*(u_\rho-\bar u_{*,\rho}),\bar p\big)\notag\\
  &=a(y_*^\dag-\bar y_{*,\rho}^\dag,R_*\bar p)
  +a\big(B_*(u_\rho-\bar u_{*,\rho}),\bar p-R_*\bar p\big)\notag\\
  &=\Big(\int_\O f R_*\bar p\,dx+\int_\Gamma u_\rho R_*\bar p\,ds\Big)-
     \Big(\int_\O f R_*\bar p\,dx+\int_\Gamma \bar u_{*,\rho} R_*\bar p\,ds\Big)\\
     &\hspace{50pt}+a\big(B_*(u_\rho-\bar u_{*,\rho}),\bar p-R_*\bar p\big)\notag\\
     &=\int_\Gamma(u_\rho-\bar u_{*,\rho})R_*\bar p\,ds
     +a\big(B_*(u_\rho-\bar u_{*,\rho}),\bar p-R_*\bar p\big).\notag
\end{align}
\par
 Let $\tilde u_\rho\in W_\rho$ be defined by
\begin{equation}\label{eq:LTwoError2}
  \tilde u_\rho=Q_\rho\bar u.
\end{equation}
 Then $\tilde u_\rho$ belongs to $U_{ad}^\rho$ by
 Remark~\ref{rem:AdmissibleSets} and we define
\begin{equation}\label{eq:tildey*}
 \tilde y_*=\tilde y_*^\dag-B_*\tilde u_\rho,
\end{equation}
where $\tilde y_{*,\rho}^\dag\in V_*$
 satisfies
\begin{equation}\label{eq:LTwoError3}
  a(\tilde y_*^\dag,z_*)=\int_\O fz_*dx
  +\int_\Gamma \tilde u_\rho z_*ds \qquad \forall\,z_*\in V_*
\end{equation}
 so that $(\tilde y_*,\tilde u_\rho)\in\bK_{*,\rho}$.
\par
 We can write
\begin{align}\label{eq:LTwoError4}
  &\|\bar y-\bar y_{*,\rho}\|_\LT^2+\gamma\|\bar u-\bar u_{*,\rho}\|_\LTG^2\notag\\
  &\hspace{40pt}=\int_\O (\bar y-\bar y_{*,\rho})(\bar y-\tilde y_*)dx+
  \gamma\int_\Gamma(\bar u-\bar u_{*,\rho})(\bar u-\tilde u_\rho)ds\\
  &\hspace{80pt}+\int_\O (\bar y-\bar y_{*,\rho})(\tilde y_*-\bar y_{*,\rho})dx+
    \gamma\int_\Gamma(\bar u-\bar u_{*,\rho})(\tilde u_\rho-\bar u_{*,\rho})ds.\notag
\end{align}
\par
 It follows from \eqref{eq:KKT}, \eqref{eq:DiscreteFirstOrderOptimality},
  \eqref{eq:LTwoError1}  and  \eqref{eq:LTwoError2}  that
\begin{align}\label{eq:LTwoError5}
  &\int_\O (\bar y-\bar y_{*,\rho})(\tilde y_*-\bar y_{*,\rho})dx+
    \gamma\int_\Gamma(\bar u-\bar u_{*,\rho})(\tilde u_\rho-\bar u_{*,\rho})ds\notag\\
    &\hspace{40pt}=\int_\O \bar y (\tilde y_*-\bar y_{*,\rho})dx
    +\gamma\int_\Gamma \bar u(\tilde u_\rho-\bar u_{*,\rho})ds\notag\\
    &\hspace{80pt}
    -\Big[\int_\O \bar y_{*,\rho}(\tilde y_*-\bar y_{*,\rho})dx
    +\gamma\int_\Gamma \bar u_{*,\rho}(\tilde u_\rho-\bar u_{*,\rho})ds\Big]\\
    &\hspace{40pt}\leq\int_\O (\bar y-\bar y_d)(\tilde y_*-\bar y_{*,\rho})dx+
    \gamma\int_\Gamma \bar u(\tilde u_\rho-\bar u_{*,\rho})ds\notag\\
    &\hspace{40pt}=\int_\Gamma(R_*\bar p+\gamma\bar u) (\tilde u_\rho-\bar u_{*,\rho})ds
    +a\big(B_*(\tu_\rho-\bar u_{*,\rho}),\bar p-R_*\bar p\big)\notag\\
    &\hspace{40pt}=\int_\Gamma (\lambda_1+\lambda_2)(\tilde u_\rho-\bar u_{*,\rho})ds
     +\int_\Gamma (R_*\bar p-\bar p)(\tilde u_\rho-\bar u_{*,\rho})ds\notag\\
     &\hspace{80pt}
     +a\big(B_*(Q_\rho\bar u-\bar u_{*,\rho}),\bar p-R_*\bar p\big).\notag
\end{align}
\par
 We have
\begin{equation}\label{eq:LTwoError6}
  \int_\Gamma (R_*\bar p-\bar p)(\tilde u_\rho-\bar u_{*,\rho})ds
  \leq \|R_*\bar p-\bar p\|_\LTG\big(\|Q_\rho \bar u-\bar u\|_\LTG
  +\|\bar u-\bar u_{*,\rho}\|_\LTG\big)
\end{equation}
 by \eqref{eq:LTwoError2}, the Cauchy-Schwarz inequality and the triangle inequality,
 and also
\begin{align}\label{eq:LTwoExtra}
  a\big(B_*(Q_\rho\bar u-\bar u_{*,\rho}),\bar p-R_*\bar p\big)
   \leq
   \|B_*\|\big(\|\bar u-Q_\rho\bar u\|_\LTG+\|\bar u-\bar u_{*,\rho}\|_\LTG\big)\|\bar p-R_*\bar p\|_a.
\end{align}
\par
 We can estimate the first term on the right-hand side of \eqref{eq:LTwoError5} by
 \eqref{eq:Signs},
 \eqref{eq:Complemntarity}, \eqref{eq:DiscreteControlConstraint}, \eqref{eq:LTwoError2}
 and the Cauchy-Schwarz inequality as follows:
\begin{align}\label{eq:LTwoError7}
  &\int_\G(\lambda_1+\lambda_2)(\tu_\rho-\bar u_{*,\rho})ds\notag\\
  &\hspace{20pt}=\int_\G \lambda_1 (\tu_\rho-\bar u)ds+\int_\G \lambda_2(\tu_\rho-\bar u)ds\notag\\
  &\hspace{50pt}+\int_\G \lambda_1(\bar u-\phi_1)ds+\int_\G\lambda_2(\bar u-\phi_2)ds\notag\\
  &\hspace{80pt}+\int_\G\lambda_1(\phi_1-Q_\rho\phi_1)ds+\int_\G\lambda_2(\phi_2-Q_\rho\phi_2)ds\notag\\
  &\hspace{110pt}+\int_\G \lambda_1(Q_\rho\phi_1-\bar u_{*,\rho})ds+
  \int_\G \lambda_2(Q_\rho\phi_2-\bar u_{*,\rho})ds\notag\\
  &\hspace{20pt}\leq \int_\G \lambda_1 (\tu_\rho-\bar u)ds+\int_\G \lambda_2(\tu_\rho-\bar u)ds\\
  &\hspace{50pt}+\int_\G\lambda_1(\phi_1-Q_\rho\phi_1)ds+\int_\G\lambda_2(\phi_2-Q_\rho\phi_2)ds\notag\\
  &\hspace{20pt}=\int_\G(\lambda_1-Q_\rho\lambda_1)(Q_\rho\bar u-\bar u)ds
  +\int_\G(\lambda_2-Q_\rho\lambda_2)(Q_\rho\bar u-\bar u)ds\notag\\
  &\hspace{50pt}+\int_\G (\lambda_1-Q_\rho\lambda_1)(\phi_1-Q_\rho\phi_1)ds+
  \int_\G (\lambda_2-Q_\rho\lambda_2)(\phi_2-Q_\rho\phi_2)ds\notag\\
  &\hspace{20pt}\leq \big(\|\lambda_1-Q_\rho\lambda_1\|_\LTG+\|\lambda_2-Q_\rho\lambda_2\|_\LTG\big)
  \|\bar u-Q_\rho\bar u\|_\LTG \notag\\
  &\hspace{50pt}+\|\lambda_1-Q_\rho\lambda_1\|_\LTG\|\phi_1-Q_\rho\phi_1\|_\LTG\notag\\
  &\hspace{80pt}
  +\|\lambda_2-Q_\rho\lambda_2\|_\LTG\|\phi_2-Q_\rho\phi_2\|_\LTG. \notag
\end{align}
\par
 Note that \eqref{eq:PDEConstraint}, \eqref{eq:Ritz},
  \eqref{eq:LTwoError2} and \eqref{eq:LTwoError3} imply
\begin{align*}
  a(R_*\bar y-\ty_*^\dag,z_*)&=a(\bar y-\ty_*^\dag,z_*)\\
  &=\Big(\int_\O fz_*dx+\int_\G \bar u z_*ds\Big)
  -\Big(\int_\O fz_*dx+\int_\G \tilde u_\rho z_*ds\Big)\\
  &=\int_\G (\bar u-Q_\rho\bar u)z_*ds \qquad\forall\,z_*\in V_*,
\end{align*}
 and hence
\begin{equation*}
  \|R_*\bar y-\ty_*^\dag\|_\LT \leq \CTR\CPF\max(\alpha^{-1},1)\|\bar u-Q_\rho\bar u\|_\LTG
\end{equation*}
 by Lemma~\ref{lem:SimpleBdd} so that
\begin{align}\label{eq:LTwoError8}
 \|\bar y-\ty_*\|_\LT&=\|\bar y-(\ty_*^\dag-B_*\tilde u_\rho)\|_\LT\notag\\
 &\leq \|\bar y-(R_*\bar y-B_*\bar u)\|_\LT+
   \|R_*\bar y-\ty_*^\dag\|_\LT+\|B_*(\bar u-\tilde u_\rho)\|_\LT\\
 &\leq \|\bar y-(R_*\bar y-B_*\bar u)\|_\LT
  +\CTR\CPF\max(\alpha^{-1},1)\|\bar u-Q_\rho\bar u\|_\LTG\notag\\
  &\hspace{50pt}+\|B_*\|(\COM)^\frac12\|\bar u-Q_\rho\bar u\|_\LTG,
 \notag
\end{align}
 where we have used  \eqref{eq:Coercive} and  \eqref{eq:LTwoError2}.
\par
 Therefore we have, by \eqref{eq:LTwoError2}, \eqref{eq:LTwoError8}
 and the Cauchy-Schwarz inequality,
\begin{align}\label{eq:LTwoError9}
  &\,\int_\O (\bar y-\bar y_{*,\rho})(\bar y-\tilde y_*)dx+
  \gamma\int_\Gamma(\bar u-\bar u_{*,\rho})(\bar u-\tilde u_\rho)ds\notag\\
  \leq&\, \|\bar y-\bar y_{*,\rho}\|_\LT
  (\CTR\COM+(\COM)^\frac12\|B_*\|)\|\bar u-Q_\rho\bar u\|_\LTG\\
  &\hspace{20pt}
  +\|\bar y-\bar y_{*,\rho}\|_\LT\|\bar y-(R_*\bar y-B_*\bar u)\|_\LT
  +\gamma\|\bar u-\bar u_{*,\rho}\|_\LTG\|\bar u-Q_\rho\bar u\|_\LTG.
  \notag
\end{align}
\par
 Putting \eqref{eq:LTwoError4}--\eqref{eq:LTwoError7} and \eqref{eq:LTwoError9}
 together, we find
\begin{align*}
   &\,\|\bar y-\bar y_{*,\rho}\|_\LT^2+\gamma\|\bar u-\bar u_{*,\rho}\|_\LTG^2\\
   \leq&\, \|R_*\bar p-\bar p\|_\LTG\big(\|Q_\rho \bar u-\bar u\|_\LTG
   +\|\bar u-\bar u_{*,\rho}\|_\LTG\big)\\
   &\hspace{10pt}+
     \|B_*\|\big(\|\bar u-Q_\rho\bar u\|_\LTG
     +\|\bar u-\bar u_{*,\rho}\|_\LTG\big)\|\bar p-R_*\bar p\|_a\\
      &\hspace{15pt}+\|\bar y-\bar y_{*,\rho}\|_\LT
   (\CTR\COM+(\COM)^\frac12\|B_*\|)\|\bar u-Q_\rho\bar u\|_\LTG\\
     &\hspace{20pt}+ \|\bar y-\bar y_{*,\rho}\|_\LT\|\bar y-(R_*\bar y-B_*\bar u)\|_\LT
  +\gamma\|\bar u-\bar u_{*,\rho}\|_\LTG\|\bar u-Q_\rho\bar u\|_\LTG\\
   &\hspace{25pt}+\big(\|\lambda_1-Q_\rho\lambda_1\|_\LTG
   +\|\lambda_2-Q_\rho\lambda_2\|_\LTG\big)
  \|\bar u-Q_\rho\bar u\|_\LTG \notag\\
  &\hspace{30pt}+\|\lambda_1-Q_\rho\lambda_1\|_\LTG\|\phi_1-Q_\rho\phi_1\|_\LTG
   +\|\lambda_2-Q_\rho\lambda_2\|_\LTG\|\phi_2-Q_\rho\phi_2\|_\LTG
\end{align*}
 which together with the inequality of arithmetic and geometric means implies
\begin{align*}
  &\,\|\bar y-\bar y_{*,\rho}\|_\LT^2+\gamma\|\bar u-\bar u_{*,\rho}\|_\LTG^2\notag\\
  \leq&\, C_\diamond\Big[\|\bar y-(R_*\bar y-B_*\bar u)\|_\LT^2
  +\gamma^{-1}\|\bar p-R_*\bar p\|_\LTG^2
   +\gamma^{-1}\|B_*\|^2\|\bar p-R_*\bar p\|_a^2\notag\\
  &\hspace{20pt}+\big(1
  +[\CTR\COM]^2+[\COM+1]\|B_*\|^2\big)\|\bar u-Q_\rho\bar u\|_\LTG^2\notag\\
  &\hspace{40pt}+\|\lambda_1-Q_\rho\lambda_1\|_\LTG^2+\|\lambda_2-Q_\rho\lambda_2\|_\LTG^2\\
  &\hspace{60pt}+\|\phi_1-Q_\rho\phi_1\|_\LTG^2+\|\phi_2-Q_\rho\phi_2\|_\LTG^2\Big],\notag
\end{align*}
 where $C_\diamond$ is a universal positive constant.
\end{proof}
\par

\section{Application to Problems with Smooth Coefficients}\label{sec:Smooth}
 In the case of smooth $A(x)$ and $\kappa(x)$, we have more information on $\bar y$,
  $\bar p$ and $\bar u$.
 For simplicity we only consider the case where $\O$ is  convex.

\subsection{Regularity Results for the Continuous Problem}\label{eq:subsubsec:Continuous}
 From the elliptic regularity result in
 \cite[Theorem~3.2.1.3]{Grisvard:1985:EPN} for homogeneous
 Neumann boundary value problems,
 we know from
 \eqref{eq:AdjointState} that
\begin{equation}\label{eq:AdjointRegularity}
  \|\bar p\|_{H^2(\O)}\leq C_\O \|\bar y-y_d\|_\LT.
\end{equation}
\par
 Let  $\Sigma$ be the set of the sides (edges for $d=2$ and faces for $d=3$)  of $\G$.
 It follows from the trace theorem (cf. \cite{ADAMS:2003:Sobolev}) that
\begin{equation}\label{eq:AdjointStatePiecewiseHOne}
  \bar p\big|_\sigma\in H^{3/2}(\sigma) \qquad\forall\,\sigma\in\Sigma.
\end{equation}
\par
 We assume that the control constraints
\begin{equation}\label{eq:PhiPiecewiseHOne}
  \text{$\phi_i$ ($i=1,2$) belong to $H^1(\sigma)$
  for every $\sigma\in\Sigma$.}\hspace{170pt}
\end{equation}
 We can then obtain from \eqref{eq:ControlAndAdjointState},
 \eqref{eq:AdjointStatePiecewiseHOne}, \eqref{eq:PhiPiecewiseHOne}
 and \cite[Lemma~7.6]{GT:2001:EllipticPDE}
  that
\begin{equation}\label{eq:ControlHOneBound}
  |\bar u|_{H^1(\sigma)}\leq \max\big(|\phi_1|_{H^1(\sigma)},|\phi_2|_{H^1(\sigma)},
  \gamma^{-1}|\bar p|_{H^1(\sigma)}\big)
  \qquad\forall\,\sigma\in\Sigma.
\end{equation}
\par
 Similarly, we have
\begin{alignat}{3}
 |\lambda_1|_{H^1(\sigma)}&\leq \max\big(|\bar p|_{H^1(\sigma)},|\phi_1|_{H^1(\sigma)}\big)&\qquad&
 \forall\,\sigma\in\Sigma,\label{eq:Lambda1Regularity}\\
  |\lambda_2|_{H^1(\sigma)}&\leq \max\big(|\bar p|_{H^1(\sigma)},|\phi_2|_{H^1(\sigma)}\big)&\qquad&
 \forall\,\sigma\in\Sigma,\label{eq:Lambda2Regularity}
\end{alignat}
 by \eqref{eq:Lambda1}, \eqref{eq:Lambda2} and the assumption that $\gamma\leq 1$.
\par
 Now we turn to the regularity of $\bar y$, which is the solution of the
 nonhomogeneous Neumann boundary value problem \eqref{eq:PDEConstraint} (with $u=\bar u$).
 It follows from \eqref{eq:ControlHOneBound} that
\begin{equation}\label{eq:StateRegularity}
  \|\bar y\|_{H^2(\O)}\leq C_\Omega \big(\|f\|_\LT
  +\sum_{\sigma\in\Sigma}\|\bar u\|_{H^{1/2}(\sigma)}\big).
\end{equation}
 For $d=2$ this regularity result follows from Theorem~1.4.3,
 Theorem~2.4.3 and Remark~2.4.5
 in \cite{Grisvard:1992:Singularities}.  For $d=3$ it follows from
 Theorem~8.1.7 in \cite{MR:2010:Polyhedral}
 and the eigenvalue estimate in \cite{Escobar:1990:Eigenvalue}, and it
 also follows from the results in
 \cite{BS:1987:Neumann}.
 %
\subsection{Results for the $P_1$ Finite Element Method} \label{subsubsec:FEM}
 We take $V_*=V_h$ to be the standard $P_1$  finite element space associated with
 a regular triangulation $\cT_h$ of $\O$ with mesh size $h$.
   We then take $\cT_\rho$ to be a refinement of the triangulation on $\G$
   induced by $\cT_h$ and denote
   the solution for the optimal control problem  by
   $(\bar y_{h,\rho},\bar u_{h,\rho}, \bar p_{h,\rho})$.
   The Ritz projection operator $R_*:\HOne\longrightarrow V_*$ is now
   denoted by $R_h:\HOne\longrightarrow V_h$,
   i.e.,
\begin{equation}\label{eq:Rh}
  a(R_h\zeta,v_h)=a(\zeta,v_h)\qquad\forall\,\zeta\in\HOne, v_h\in V_h.
\end{equation}
 In this case we take $B_*=0$.
\par
 It follows from \eqref{eq:AdjointRegularity},
 \eqref{eq:ControlHOneBound}--\eqref{eq:Lambda2Regularity},
 \eqref{eq:StateRegularity}
  and
 standard finite element estimates (cf. \cite{Ciarlet:1978:FEM,BScott:2008:FEM}) that
\begin{align}
  \|\bar y-R_h\bar y\|_\LT&\leq Ch^{2}\|\bar y\|_{H^{2}(\O)},
  \label{eq:SmoothStatsLTwoError}\\
   \|\bar p-R_h\bar p\|_\HOne&\leq Ch\|\bar p\|_{H^{2}(\O)},
  \label{eq:SmoothAdjointStateEnergyError}\\
   \|\bar p-R_h\bar p\|_\LT&\leq Ch^{2}\|\bar p\|_{H^{2}(\O)},
  \label{eq:SmoothAdjointStatsLTwoError}\\
  \|\bar u-Q_\rho\bar u\|_\LTG&\leq C \rho\sum_{\sigma\in\Sigma}|\bar u|_{H^1(\sigma)},
  \label{eq:SmoothControlLTwoError}\\
  \|\lambda_1-Q_\rho\lambda_1\|_\LTG&\leq C\rho\sum_{\sigma\in\Sigma}
  |\lambda_1|_{H^1(\sigma)},
  \label{eq:SmoothLambda1LTwoError}\\
   \|\lambda_2-Q_\rho\lambda_2\|_\LTG&\leq C\rho\sum_{\sigma\in\Sigma}
   |\lambda_2|_{H^1(\sigma)},
  \label{eq:SmoothLambda2LTwoError}\\
  \|\phi_1-Q_\rho\phi_1\|_\LTG&\leq C\rho\sum_{\sigma\in\Sigma}|\phi_1|_{H^1(\sigma)},
  \label{eq:Smoothphi1LTwoError}\\
    \|\phi_2-Q_\rho\phi_2\|_\LTG&\leq C\rho\sum_{\sigma\in\Sigma}|\phi_2|_{H^1(\sigma)},
  \label{eq:Smoothphi2LTwoError}
\end{align}
 where the positive constant $C$ only depends on the shape regularity of $\cT_h$ and $\cT_\rho$.
\par
 Note that \eqref{eq:SmoothAdjointStateEnergyError} and \eqref{eq:SmoothAdjointStatsLTwoError}
 imply
\begin{equation}\label{eq:SmoothAdjointBoundaryLTwoError}
  \|\bar p-R_h\bar p\|_\LTG\leq Ch^{3/2}\|\bar p\|_{H^2(\O)}
\end{equation}
through the trace inequality (cf. \cite[Theorem~1.6.6]{BScott:2008:FEM})
\begin{equation*}
  \|v\|_\LTG\leq C_\O \|v\|_\LT^\frac12|v|_\HOne^\frac12 \qquad\forall\,v\in \HOne.
\end{equation*}
\par
 Combining Theorem~\ref{thm:LTwoError},
 \eqref{eq:ControlHOneBound}--\eqref{eq:Lambda2Regularity},
 \eqref{eq:SmoothStatsLTwoError},
 \eqref{eq:SmoothControlLTwoError}--\eqref{eq:SmoothAdjointBoundaryLTwoError}
  and Remark~\ref{rem:Bounds}, we conclude
 that
\begin{align}
  \|\bar y-\bar y_{h,\rho}\|_\LT+\|\bar u-\bar u_{h,\rho}\|_\LTG
    &\leq C_\ddag(h^{3/2}+\rho),\label{eq:SmoothLTwoError}
\end{align}
 where the positive constant $C_\ddag$ only depends on $\alpha^{-1}$, $\gamma^{-1}$,
  $\|f\|_\LT$,
 $\|y_d\|_\LT$, $\sum_{\sigma\in\Sigma}\|\phi_i\|_{H^1(\sigma)}$
 ($i=1,2$),
  and the shape regularity of $\cT_h$ and $\cT_\rho$.
\par
 In the case where $\cT_\rho$ is the triangulation on $\G$ induced by $\cT_h$, the $L_2$
 error estimate becomes
\begin{equation*}
    \|y-\bar y_{h,\rho}\|_\LT+\|\bar u-\bar u_{h,\rho}\|_\LTG
  \leq C h,
\end{equation*}
 and we have recovered the two dimensional result in \cite{CMT:2005:Neumann}
  for the model problem.
%
\subsection{Variational Discretization Concept}\label{subsubsec:Hinze}
 If  we adopt the variational discretization concept of Hinze
  (cf. \cite[Chapter~3]{HPUU:2009:Book}) and do not discretize the control variable, then
   we can replace $W_\rho$  by $\LTG$ and
 $Q_\rho:\LTG\longrightarrow W_\rho$  by the identity operator on $\LTG$.
\par
 Let the solution of the optimal control problem be denoted by
 $(\bar y_h,\bar u_h,\bar p_h)$.
 According to Theorem~\ref{thm:LTwoError},
 the estimate for the $L_2$ errors becomes
\begin{equation*}
  \|\bar y-\bar y_h\|_\LT+\|\bar u-\bar u_h\|_\LTG\leq
   C_\flat\big(\|\bar y-R_h\bar y\|_\LT
  +\|\bar p-R_h\bar p\|_\LTG\big)
\end{equation*}
 which implies the estimate in \cite[display (10)]{HM:2009:Neumann}.
\par
 The estimate corresponding to
  \eqref{eq:SmoothLTwoError}  now reads
\begin{equation*}
  \|\bar y-\bar y_{h,\rho}\|_\LT+\|\bar u-\bar u_{h,\rho}\|_\LTG
  \leq C h^{3/2}
\end{equation*}
 and we have recovered the result in \cite{CM:2008:Neumann}
 for the model problem in two dimensions.
%
\section{Application to Problems with Rough Coefficients}\label{sec:Rough}
 We consider the optimal control problem
 \eqref{eq:LQ}--\eqref{eq:ControlConstraint} on
 an arbitrary polygonal (resp., polyhedral) domain in $\R^2$ (resp., $\R^3$).
 In the case of rough coefficients, where \eqref{eq:RoughCoefficient}
 is the only assumption on
 the matrix $A(x)$, we have $\bar y,\bar p\in\HOne$ and nothing more.
 Therefore we can only
 assert that the trace of $\bar p$ on $\G$ belongs to $H^{1/2}(\G)$.
\par
 Under the assumption that the control constraints
\begin{equation}\label{eq:RoughConstraints}
  \text{$\phi_1$ and $\phi_2$ belong to $\HH$,}
\end{equation}
 we have the following regularity estimates from
 \eqref{eq:ControlAndAdjointState}--\eqref{eq:Lambda2}:
\begin{align}
  \|\bar u\|_\HH&\leq C_\O\max\big(\|\phi_1\|_\HH+\|\phi_2\|_\HH
  +\gamma^{-1}\|\bar p\|_\HH\big),\label{eq:RoughControlRegularity}\\
  \|\lambda_1\|_\HH&\leq C_\O\max\big(\|\bar p\|_\HH+\|\phi_1\|_\HH\big),
  \label{eq:RoughLambda1Regularity}\\
  \|\lambda_2\|_\HH&\leq C_\O\max\big(\|\bar p\|_\HH+\|\phi_2\|_\HH\big).
  \label{eq:RoughLambda2Regularity}
\end{align}
\par
 It follows from
 \eqref{eq:RoughConstraints}--\eqref{eq:RoughLambda2Regularity}
 and standard finite element estimates that
\begin{align}
  \|\phi_1-Q_\rho\phi_1\|_\LTG&\leq C\rho^\frac12\|\phi_1\|_\HH,\label{eq:Roughphi1LTwoError}\\
  \|\phi_2-Q_\rho\phi_2\|_\LTG&\leq C\rho^\frac12\|\phi_2\|_\HH,\label{eq:Roughphi2LTwoError}\\
  \|\bar u-Q_\rho\bar u\|_\LTG&\leq C\rho^\frac12\|\bar u\|_\HH,\label{eq:RoughControlLTwoError}\\
  \|\lambda_1-Q_\rho\lambda_1\|_\LTG&\leq C\rho^\frac12\|\lambda_1\|_\HH,\label{eq:RoughLambda1LTwoError}\\
  \|\lambda_2-Q_\rho\lambda_2\|_\LTG&\leq C\rho^\frac12\|\lambda_2\|_\HH,\label{eq:RoughLambda2LTwoError}
\end{align}
 where the positive constant $C$ only depends on the shape regularity of $\cT_\rho$.
%
\subsection{A Standard $P_1$ Finite Element Method}\label{subsubsec:StandardFEM}
 If we solve the optimal control problem by taking $V_*$ to be a standard $P_1$
 finite element space $V_h$ associated with a triangulation
 $\cT_h$ of $\O$ and denote the solution of the discrete problem by
 $(\bar y_{h,\rho},\bar u_{h,\rho},\bar p_{h,\rho})$, we have
\begin{equation}\label{eq:RoughStandardLTwoError}
  \|\bar y-\bar y_{h,\rho}\|_\LT+\|\bar u-\bar u_{h,\rho}\|_\LT
  \leq C_\spadesuit\big(\|\bar y-R_h\bar y\|_\LT+\|\bar p-R_h\bar p\|_\LTG+\rho^\frac12\big)
\end{equation}
 by  Theorem~\ref{thm:LTwoError} and
 \eqref{eq:Roughphi1LTwoError}--\eqref{eq:RoughLambda2LTwoError},
  where, in view of
  \eqref{eq:RoughControlRegularity}--\eqref{eq:RoughLambda2Regularity}
  and Remark~\ref{rem:Bounds}, the positive constant $C_\spadesuit$ only depends on
  $\alpha^{-1}$, $\gamma^{-1}$, $\|y_d\|_\LT$, $\|f\|_\LT$, $\|\phi_1\|_{H^{1/2}(\G)}$,
  $\|\phi_2\|_{H^{1/2}(\G)}$
  and the shape regularity of $\cT_\rho$.
\par
 Since we only have $\bar p\in \HOne$ and nothing more, we may as well
 replace $\|\bar p-R_h\bar p\|_\LTG$
 by $\CTR\|\bar p-R_h\bar p\|_\HOne\leq \CTR\sqrt{\COM}\|\bar p-R_h\bar p\|_a
 $ so that \eqref{eq:RoughStandardLTwoError} is replaced by
\begin{equation}\label{eq:NewRoughStandardLTwoError}
  \|\bar y-\bar y_{h,\rho}\|_\LT+\|\bar u-\bar u_{h,\rho}\|_\LT
  \leq C_\spadesuit\big(\|\bar y-R_h\bar y\|_\LT+\|\bar p-R_h\bar p\|_a+\rho^\frac12\big).
\end{equation}
\par
 We can also bound the energy norm error of the discrete adjoint state as follows.
 We have
\begin{equation}\label{eq:AdjointStateError1}
  \|\bar p-\bar p_{h,\rho}\|_a \leq \|\bar p-R_h\bar p\|_a+
  \|R_h\bar p-\bar p_{h,\rho}\|_a,
\end{equation}
 and
\begin{equation*}
  a(q_h,R_h\bar p-\bar p_{h,\rho})=a(q_h,\bar p-\bar p_{h,\rho})
                   =\int_\O (\bar y-\bar y_{h,\rho})q_h  \qquad\forall\,q_h\in V_h
\end{equation*}
 by \eqref{eq:AdjointState}, \eqref{eq:DiscreteAdjointState} and \eqref{eq:Rh} so that
\begin{equation}\label{eq:AdjointStateError2}
  \|R_h\bar p-\bar p_{h,\rho}\|_a\leq \sqrt{\COM}\|\bar y-\bar y_{h,\rho}\|_\LT
\end{equation}
 by Lemma~\ref{lem:SimpleBdd}.
\par
 It follows from \eqref{eq:NewRoughStandardLTwoError}--\eqref{eq:AdjointStateError2} that
\begin{equation}\label{eq:FinalRCStandardBound}
  \|\bar y-\bar y_{h,\rho}\|_\LT+\|\bar u-\bar u_{h,\rho}\|_\LT+\|\bar p-\bar p_{h,\rho}\|_a
  \leq C_\clubsuit \big(\|\bar y-R_h\bar y\|_\LT+\|\bar p-R_h\bar p\|_a+\rho^\frac12\big),
\end{equation}
where the positive constant $C_\clubsuit$ only depends on
  $\alpha^{-1}$, $\gamma^{-1}$, $\|y_d\|_\LT$, $\|f\|_\LT$, $\|\phi_1\|_{H^{1/2}(\G)}$,
  $\|\phi_2\|_{H^{1/2}(\G)}$
  and the shape regularity of $\cT_\rho$.
\par
 By a density argument, we have
\begin{equation*}
  \lim_{h\downarrow0}\|\zeta-R_h\zeta\|_\HOne=0 \qquad\text{for any $\zeta\in \HOne$},
\end{equation*}
 and hence
\begin{equation*}
  \lim_{h,\rho\downarrow 0}\big(\|\bar y-\bar y_{h,\rho}\|_\LT
  +\|\bar u-\bar u_{h,\rho}\|_\LTG
  +\|\bar p-R_h\bar p\|_a\big)=0
\end{equation*}
 by \eqref{eq:FinalRCStandardBound}.
 However the convergence in $h$ can be arbitrarily slow (cf. \cite{BO:2000:Bad})
 so that a good
 approximate solution may require a very small $h$.  This is computationally expensive,
 especially if the optimal control problem \eqref{eq:LQ}--\eqref{eq:ControlConstraint}
  has to be solved repeatedly for different
 $f$. Therefore it is desirable to
  improve the performance of the
 finite
 element method on coarse meshes by taking
 $V_*$ to be a multiscale finite element space.
\par
 The  result below is a  sort of converse to \eqref{eq:FinalRCStandardBound}.
   It is useful for  the analysis of the
  multiscale finite element methods in Section~\ref{subsubsec:MSFEM}.
\begin{lemma}\label{lem:Converse}
    There exists a positive constant $C_\natural$ depending only on $\alpha^{-1}$ such that
\begin{equation}\label{eq:Converse}
  \|\bar y-R_h\bar y\|_\LT+\|\bar p-R_h\bar p\|_a\leq C_\natural
  \big(\|\bar y-\bar y_{h,\rho}\|_\LT+\|\bar u-\bar u_{h,\rho}\|_\LTG
  +\|\bar p-\bar p_{h,\rho}\|_a\big).
\end{equation}
\end{lemma}
\begin{proof}
  We begin with a triangle inequality
\begin{align}\label{eq:Converse1}
  &\|\bar y-R_h\bar y\|_\LT+\|\bar p-R_h\bar p\|_a\\
  &\hspace{40pt}
   \leq\|\bar y-\bar y_{h,\rho}\|_\LT+\|\bar y_{h,\rho}-R_h\bar y\|_\LT+
    \|\bar p-\bar p_{h,\rho}\|_a+\|\bar p_{h,\rho}-R_h\bar p\|_a.
   \notag
\end{align}
\par
 Observe that
\begin{equation*}
  a(\bar y_{h,\rho}-R_h\bar y,z_h)=a(\bar y_{h,\rho}-\bar y,z_h)=\int_\G (\bar u_{h,\rho}-\bar u)z_hds
  \qquad\forall\,z_h\in V_h
\end{equation*}
 by \eqref{eq:PDEConstraint}, \eqref{eq:DiscretePDEConstraint}
 (where $\bar y_{h,\rho}=\bar y_{h,\rho}^\dag$ because $B_*=0$) and \eqref{eq:Rh},
 and hence
\begin{equation}\label{eq:Converse2}
  \|\bar y_{h,\rho}-R_h\bar y\|_\HOne\leq \CTR\COM\|\bar u-\bar u_{h,\rho}\|_\LTG
\end{equation}
 by Lemma~\ref{lem:SimpleBdd}.
\par
 The estimate \eqref{eq:Converse} follows from \eqref{eq:AdjointStateError2},
 \eqref{eq:Converse1} and \eqref{eq:Converse2}.
\end{proof}
%
\subsection{Multiscale Finite Element Methods}\label{subsubsec:MSFEM}
 Let $\cT_H$ be a triangulation of $\O$ with mesh size $H$ so that $\cT_h$ ($h\ll H$)
  is a refinement of $\cT_H$.  The $P_1$ finite element space associated with
   $\cT_H$ (resp., $\cT_h$) is denoted by
   $V_H$ (resp., $V_h$).
\par
 Let $V_*=\MS\subset V_h$ be the ideal multiscale finite element
  space in \cite{HM:2014:LODBdry} obtained by the orthogonal decomposition methodology.
    The construction of $\MS$ involves a quasi-local interpolation operator
  $\IH:\HOne\longrightarrow V_H$ that satisfies
\begin{equation}\label{eq:IHBdd}
  H_T^{-1}\|v-\IH v\|_{L_2(T)}+\|\nabla(v-\IH v)\|_{L_2(T)}\leq C_{\IH}\|\nabla v\|_{L_2(\omega_T)}
  \qquad\forall\,T\in\cT_H,
\end{equation}
 where $H_T$ is the diameter of $T$, $\omega_T$ is the star of
 $T$ and the positive constant $C_{\IH}$
 only depends on the shape regularity of $\cT_H$.
 \par
 Let $W_h=\{w_h\in V_h:\,\IH w_h=0\}$ be the kernel of $\IH$ on $V_h$.
 The space $\MS\subset V_h$ is defined to be the
 orthogonal complement of $W_h$ with respect to the bilinear
 form $a(\cdot,\cdot)$, i.e.,
\begin{equation}\label{eq:IdealMS}
  \MS=\{v_h\in V_h:\,a(v_h,w_h)=0 \quad\forall\,w_h\in W_h\}.
\end{equation}
\begin{remark}\label{rem:Ideal}\rm
  The dimension of $\MS$ is identical to the
 dimension of $V_H$.  The construction of the basis
 of $\MS$, which involves solving problems associated with $\cT_h$,
  can be carried out off-line
  and the
 on-line computation involving the space $\MS$ with a small dimension is fast.
\end{remark}
\par
 The operator $B_*:\LTG\longrightarrow W_h$ is defined by
\begin{equation}\label{eq:B*Def}
  a(B_*q,w_h)=-\int_\G q w_h ds\qquad\forall\,q\in\LTG,\;w_h\in W_h.
\end{equation}
 It follows from \eqref{eq:IdealMS} that \eqref{eq:Orhtogonality} is satisfied.
\par
 We deduce from Lemma~\ref{lem:SimpleBdd} that
\begin{equation*}
  \|B_*q\|_a\leq \sqrt{\COM}\,\CTR\|q\|_\LTG \qquad\forall\,q\in\LTG
\end{equation*}
 so that
\begin{equation}\label{eq:BStarNorm}
  \|B_*\|\leq \CTR\sqrt{\COM}.
\end{equation}

\par
 The following is the salient property of $\MS$ that follows from
 \eqref{eq:IHBdd}--\eqref{eq:B*Def}
 (cf. \cite[Section~3.3]{HM:2014:LODBdry}).
\begin{lemma}\label{lem:HM2014}
 Let $v_h\in V_h$ and $\MSv\in\MS$ satisfy
\begin{alignat*}{3}
  a(v_h,w_h)&=\int_\O Fw_hdx+\int_\G Gw_hds&\qquad&\forall\,w_h\in V_h,\\
  a(\MSv,\MSw)&=\int_\O F\MSw dx+\int_\G G\MSw ds&\qquad&\forall\,\MSw\in\MS,
\end{alignat*}
where $F\in\LT$ and $G\in\LTG$.
 Then we have
\begin{align*}
  \|v_h-(\MSv-B_*G)\|_a&\leq C_\diamondsuit H\big(\|F\|_\LT+\|G\|_\LTG\big),\\
  \|v_h-(\MSv-B_*G)\|_\LT&\leq (C_\diamondsuit H)^2\big(\|F\|_\LT+\|G\|_\LTG\big),
\end{align*}
 where the positive constant $C_\diamondsuit$ is independent of the mesh sizes.
\end{lemma}
\par
 Let the solution of the discrete problem with $V_*=\MS$ be denoted by
 $(\bMSy,\bMSu,\bMSp)$.
\begin{theorem}\label{thm:MSError}
 There exist two positive constants $C_1$ and $C_2$ independent of
 the mesh sizes such that
\begin{align}\label{eq:MSError}
  &\|\bar y-\bMSy\|_\LT+\|\bar u-\bMSu\|_\LTG \\
  &\hspace{40pt}\leq C_1\big(\|\bar y-\bar y_{h,\rho}\|_\LT+\|\bar u-\bar u_{h,\rho}\|_\LTG
      +\|\bar p-\bar p_{h,\rho}\|_a\big)
      +C_2(H+\rho^\frac12),\notag
\end{align}
 where $(\bar y_{h,\rho},\bar u_{h,\rho},\bar p_{h,\rho})$ is the approximate solution
 from  Section~\ref{subsubsec:StandardFEM} computed by the
 $P_1$ finite element method on the fine mesh $\cT_h$.
\end{theorem}
\begin{proof}
 According to Theorem~\ref{thm:LTwoError}, \eqref{eq:Coercive}, \eqref{eq:Trace},
 \eqref{eq:Roughphi1LTwoError}--\eqref{eq:RoughLambda2LTwoError}
 and \eqref{eq:BStarNorm},
 we have
\begin{align}\label{eq:AbstractMSError}
  &\|\bar y-\bMSy\|_\LT+\|\bar u-\bMSu\|_\LTG\\
  &\hspace{40pt}\leq C_\Box\big(\|\bar y-(\RMS\bar y-B_*\bar u)\|_\LT
  +\|\bar p-\RMS\bar p\|_a+\rho^\frac12\big),\notag
\end{align}
 where  $\RMS:\HOne\longrightarrow \MS$ is the Ritz projection operator, i.e,
\begin{equation}\label{eq:MSRitz}
  a(\RMS\zeta,\MSv)=a(\zeta,\MSv) \qquad\forall\,\zeta\in\HOne,\MSw\in\MS,
\end{equation}
 and the positive constant $C_\Box$ only depends on
 $\alpha^{-1}$, $\gamma^{-1}$, $\|y_d\|_\LT$, $\|f\|_\LT$, $\|\phi_1\|_{H^{1/2}(\G)}$,
  $\|\phi_2\|_{H^{1/2}(\G)}$
  and the shape regularity of $\cT_\rho$.
\par
 We can estimate the terms on the right-hand side of \eqref{eq:AbstractMSError} as follows.
\par
 In view of \eqref{eq:PDEConstraint}, \eqref{eq:Rh}
 and \eqref{eq:MSRitz}, we have
\begin{alignat*}{3}
  a(R_h\bar y,w_h)&=a(\bar y,w_h)=\int_\O fw_hdx+\int_\G \bar u w_h ds 
  &\qquad&\forall\,w_h\in V_h,\\
  a(\RMS\bar y,\MSw)&=a(\bar y,\MSw)=\int_\O f\MSw dx+\int_\G \bar u \MSw ds
   &\qquad&\forall\,\MSw\in \MS.
\end{alignat*}
 It then follows from Lemma~\ref{lem:HM2014} that
\begin{equation*}
  \|R_h\bar y-(\RMS \bar y-B_*\bar u)\|_\LT\leq 
  (C_\diamondsuit H)^2\big(\|f\|_\LT+\|\bar u\|_\LTG\big)
\end{equation*}
 and hence
 \begin{equation}\label{eq:Comparison1}
   \|\bar y-(\RMS\bar y-B_*\bar u)\|_\LT\leq \|\bar y-R_h\bar y\|_\LT
   +(C_\diamondsuit H)^2(\|f\|_\LT+\|\bar u\|_\LTG\big).
 \end{equation}
 \par
  Similarly from \eqref{eq:AdjointState}, \eqref{eq:DiscreteAdjointState},
  \eqref{eq:Rh}
 and \eqref{eq:MSRitz}, we see that
\begin{alignat*}{3}
  a(z_h,R_h\bar p)&=a(z_h,\bar p)=\int_\O (\bar y-y_d)z_h dx &\qquad&\forall\,z_h\in V_h,\\
  a(\MSz,\RMS\bar p)&=a(\MSz,\bar p)=\int_\O (\bar y-y_d)\MSz dx &\qquad&\forall\,\MSz\in\MS,
\end{alignat*}
 and hence
\begin{equation*}
  \|R_h\bar p-\RMS\bar p\|_a\leq C_\diamondsuit H\|\bar y-y_d\|_\LT\label{eq:ASEnergyComp}
\end{equation*}
 by Lemma~\ref{lem:HM2014}.
 Consequently we have
\begin{equation}\label{eq:Comparison2}
   \|\bar p-\RMS\bar p\|_a\leq \|\bar p-R_h\bar p\|_a
   +(C_\diamondsuit H)\|\bar y-y_d\|_\LT.
\end{equation}
\par
 Combining
 \eqref{eq:AbstractMSError}, \eqref{eq:Comparison1} and \eqref{eq:Comparison2},
 we find by Lemma~\ref{lem:Converse}
\begin{align*}
  &\|\bar y-\bMSy\|_\LT+\|\bar u-\bMSu\|_\LTG\notag\\
  &\hspace{40pt}\leq C_\Box \big(\|\bar y-R_h\bar y\|_\LT+\|\bar p-R_h\bar p\|_a\big)
  +C_\heartsuit (H+\rho^\frac12)\\
  &\hspace{40pt}\leq  C_\Box C_\natural\big(\|\bar y-\bar y_{h,\rho}\|_\LT+
   \|\bar u-\bar u_{h,\rho}\|_\LTG+\|\bar p-\bar p_{h,\rho}\|_a\big)
   +C_\heartsuit (H+\rho^\frac12),
  \notag
\end{align*}
 where the positive constant $C_\heartsuit$ is independent of the mesh sizes.
\end{proof}
\par
  From \eqref{eq:MSError}
  we can say that, up to an $O(H+\rho^\frac12)$ error, the  
  solution obtained by the discretization
 based on $\MS\times W_\rho$ on a coarse mesh $\cT_H$ is 
 qualitatively similar to the error of
  the approximate solution obtained by the  standard $P_1$ 
  finite element method based
  on $V_h\times W_\rho$ for a potentially very small $h$.
 \par
  Alternatively, for a chosen tolerance tol, one can choose
 \begin{equation}\label{eq:rhoChoice}
 \rho^\frac12\leq \frac{\mathrm{tol}}{2C_\clubsuit}
 \end{equation}
  and a sufficiently small $h$ such that
 \begin{equation}\label{eq:hChoice}
  \|\bar y-R_h\bar y\|_\LT+\|\bar p-R_h\bar p\|_a\leq \frac{\mathrm{tol}}{2C_\clubsuit}.
 \end{equation}
 %
  It  follows from \eqref{eq:FinalRCStandardBound}, \eqref{eq:rhoChoice} and
  \eqref{eq:hChoice} that
 \begin{equation}\label{eq:FineMeshError}
  \|\bar y-\bar y_{h,\rho}\|_\LT+\|\bar u-\bar u_{h,\rho}\|_\LTG
  +\|\bar p-\bar p_{h,\rho}\|_a\leq\mathrm{tol}.
 \end{equation}
  We then derive from Theorem~\ref{thm:MSError},  \eqref{eq:rhoChoice}
  and  \eqref{eq:FineMeshError} that
 \begin{equation*}
 \|\bar y-\bMSy\|_\LT+\|\bar u-\bMSu\|_\LTG
 \leq\Big(C_1+\frac{C_2}{2C_\clubsuit}\Big)\mathrm{tol}+C_2H.
 \end{equation*}
  Hence we can also say that up to a multiple of the tolerance, the multiscale
  finite element method based on a coarse mesh $\cT_H$ has $O(H)$ convergence.

\vspace{1in}
   Therefore once a fine mesh has been determined that can resolve the multiscale features
  in the partial differential equation constraint \eqref{eq:PDEConstraint} to
   provide an acceptable
  approximate solution for the optimal control problem, one can switch
   to the multiscale finite element method
 on coarse meshes to obtain good approximate
   solutions  efficiently by solving  problems of much smaller dimensions.
\par
 The computation of the basis functions of $\MS$, which involves solving problems on $V_h$,
 may still be too expensive.  This can be remedied by the Local Orthogonal Decomposition
 (LOD) approach (cf.~\cite{MP:2014:LOD,MP:2021:LOD}), where 
 the computation is performed on patches around the
 vertices of $\cT_H$.  The difference between the solution computed by the LOD method
 and the solution computed by
  the ideal multiscale finite element method based on $\MS$ is $O(H^2)$ in $\|\cdot\|_\LT$
  and $O(H)$ in $\|\cdot\|_a$
   (cf. \cite[Lemma~3.6 and Conclusion~3.9]{HM:2014:LODBdry}), provided that the number of layers
  in the local patches is sufficiently large. Therefore
  Theorem~\ref{thm:MSError}   remains valid for the solution
  obtained by the LOD method in \cite{HM:2014:LODBdry}.
\section{Concluding Remarks}\label{sec:Conclusions}
 We have developed an abstract error analysis
 of finite element methods for a linear-quadratic elliptic Neumann boundary control problem
 purely in terms of
 the Ritz projection operator for $\HOne$, the $L_2$ projection operator for $\LTG$ and
 a correction operator for the  Neumann boundary data,
  which makes the analysis
 applicable to standard finite element methods when the coefficients of the elliptic operator
 are smooth and  to multiscale finite element methods when the coefficients are rough.
\par
For simplicity we have assumed that the bilinear form $a(\cdot,\cdot)$
in \eqref{eq:aDef} is symmetric.
  But the estimates in Section~\ref{sec:Analysis}
 can be extended to a
 nonsymmetric $a(\cdot,\cdot)$ by replacing the term $R_*\bar p$ with the term
 $S_*\bar p$, where $S_*:\HOne\longrightarrow V_*$ is defined by
\begin{equation*}
  a(q_*,S_*\zeta)=a(q_*,\zeta) \qquad\forall\,q_*\in V_*.
\end{equation*}
\par
 The results in this paper can also be applied to other multiscale finite element methods such as
 the one in the recent work \cite{YC:2023:Bdry}.
\par
 Finally we note that it is more challenging to  analyze multiscale finite element methods for
 elliptic  Dirichlet boundary control problems
 (cf. \cite{CR:2006:Dirichlet,DGH:2009:Dirichlet,MRV:2013:Dirichlet}) due to the
 ultra weak formulation of the state equation and the fact that
 the adjoint state only belongs to $\HOne$ under the assumption of
 rough coefficients.
\par\medskip
\centerline{Acknowledgement}
\par\medskip
 The authors would like to thank Johannes Pfefferer, J\"urgen Rossmann and Boris Vexler
 for helpful discussions on the regularity result for the nonhomogeneous Neumann boundary
 value problem in \cite{MR:2010:Polyhedral}, and to thank Martin Costabel for pointing to
 the reference \cite{BS:1987:Neumann}.
\appendix
\section{Some Estimates for $\bar y$, $\bar u$ and $\bar p$}\label{appendix:Bounds}
%
 We can take $u_1=\phi_1$ and $y_1\in\HOne$ defined by \eqref{eq:PDEConstraint} so that
 $(y_1,u_1)\in\bK$ and
\begin{equation*}
  a(y_1,z)=\int_\O fz\,dx+\int_\G \phi_1 z\,ds \qquad\forall\,z\in\HOne,
\end{equation*}
 which, in view of Lemma~\ref{lem:SimpleBdd} with $V=\HOne$, implies
\begin{equation}\label{eq:Bdd1}
  \|y_1\|_\HOne\leq \COM\big(\|f\|_\LT+\CTR\|\phi_1\|_\LTG\big).
\end{equation}
\par
 It then follows from \eqref{eq:Bdd1} and
 the relation
\begin{equation*}
  \|\bar y- y_d\|_\LT^2+\gamma\|\bar u\|_\LTG^2 \leq \|y_1-y_d\|_\LT^2+\gamma\|u_1\|_\LTG^2
  \leq \|y_1-y_d\|_\LT^2+\|\phi_1\|_\LTG^2
\end{equation*}
  that
\begin{align}\label{eq:Bdd2}
  \|\bar y-y_d\|_\LT^2+\gamma\|\bar u\|_\LTG^2
  &\leq 2\|y_d\|_\LT^2+
  4\big[\COM\big]^2\|f\|_\LT^2\\
  &\hspace{60pt}+\big(4\big[\CTR\COM\big]^2+1\big)\|\phi_1\|_\LTG^2.\notag
\end{align}
\par
 Similarly we have
\begin{align}\label{eq:Bdd3}
  \|\bar y-y_d\|_\LT^2+\gamma\|\bar u\|_\LTG^2
  &\leq 2\|y_d\|_\LT^2+
   4\big[\COM\big]^2\|f\|_\LT^2\\
  &\hspace{60pt}+\big(4\big[\CTR\COM\big]^2+1\big)\|\phi_2\|_\LTG^2.\notag
\end{align}
\par
 We conclude from \eqref{eq:Bdd2} and \eqref{eq:Bdd3} that
\begin{align}
  \|\bar y-y_d\|_\LT&\leq C_\sharp^\frac12,\label{eq:Bdd4}\\
  \|\bar u\|_\LTG&\leq \gamma^{-\frac12}C_\sharp^\frac12,\label{eq:Bdd5},
\end{align}
 where
\begin{align}\label{eq:CSharp}
  C_\sharp&=
  2\|y_d\|_\LT^2+
   4\big[\COM\big]^2\|f\|_\LT^2\\
  &\hspace{60pt}+\big(4\big[\CTR\COM\big]^2+1\big)
  \min\big(\|\phi_1\|_\LTG^2,\|\phi_2\|_\LTG^2\big).\notag
\end{align}
\par
 From \eqref{eq:PDEConstraint},  Lemma~\ref{lem:SimpleBdd} (with $V=\HOne$)
  and \eqref{eq:Bdd5}, we immediately have
\begin{equation}\label{eq:Bdd6}
  \|\bar y\|_\HOne\leq \COM \|\bar u\|_\LTG\leq \gamma^{-1}\COM C_\sharp.
\end{equation}
 Similarly,
 it follows from \eqref{eq:AdjointState},  Lemma~\ref{lem:SimpleBdd} (with $V=\HOne$)
 and \eqref{eq:Bdd4} that
\begin{equation}\label{eq:Bdd7}
  \|\bar p\|_\HOne\leq \COM\|\bar y-y_d\|_\LT \leq \COM C_\sharp^\frac12.
\end{equation}
\begin{remark}\label{rem:Bounds}\rm
  One can see from \eqref{eq:Bdd4}--\eqref{eq:Bdd7} that the quantities $\|\bar y- y_d\|_\LT$,
  $\|\bar u\|_\LTG$, $\|\bar y\|_\HOne$ and $\|\bar p\|_\HOne$ are bounded by constants that
  only depend on $\alpha^{-1}$, $\gamma^{-1}$, $\|f\|_\LT$, $\|y_d\|_\LT$, 
  $\|\phi_1\|_\LTG$ and $\|\phi_2\|_\LTG$.
\end{remark}

\end{document}